\documentclass{aa}  

\usepackage{natbib}
\bibpunct{(}{)}{;}{a}{}{,} 

\usepackage{graphicx}
\usepackage{tikz}
\usetikzlibrary{arrows,shapes,automata,backgrounds,petri,positioning}

\usepackage{amssymb,amsmath,bm}
\usepackage{delarray}
\usepackage{mathrsfs}
\usepackage[cal=boondox]{mathalfa}
\usepackage{xcolor}
\usepackage{soul}
\usepackage[varg]{txfonts}
\usepackage{hyperref}
\hypersetup{
  colorlinks   = true, 
  urlcolor     = blue, 
  linkcolor    = blue, 
  citecolor   = blue 
}

\usepackage{placeins}

\usepackage[ruled,vlined]{algorithm2e}

\SetCommentSty{mycommfont}

\newcommand{\rr}{{\mathbf r}}
\newcommand{\xx}{{\mathbf x}}
\newcommand{\bb}{{\mathbf b}}
\newcommand{\vv}{{\mathbf v}}
\newcommand{\ww}{{\mathbf w}}
\newcommand{\zz}{{\mathbf z}}

\begin{document} 



\title{Numerical solutions to linear transfer problems\\ 
of polarized radiation}
\subtitle{II. Krylov methods and matrix-free implementation}

\author{Pietro Benedusi\inst{1}
        \and
       Gioele Janett\inst{2,1}
        \and
       Luca Belluzzi\inst{2,3,1}
       \and
       Rolf Krause\inst{1}}
\institute{
Euler Institute, Universit\`a della Svizzera italiana (USI), CH-6900 Lugano, Switzerland
\and
Istituto Ricerche Solari (IRSOL), Università della Svizzera italiana (USI), CH-6605 Locarno-Monti, Switzerland
\and
Leibniz-Institut f\"ur Sonnenphysik (KIS), D-79104 Freiburg i.~Br.,Germany
\\      \email{pietro.benedusi@usi.ch}
}
 
\abstract
{
Numerical solutions to transfer problems of polarized radiation in solar and stellar atmospheres commonly rely on stationary iterative methods, which often perform poorly when applied to large problems.
In recent times, stationary iterative methods have been replaced by state-of-the-art preconditioned Krylov iterative methods for many applications. However, a general description and a convergence analysis of Krylov methods in the polarized radiative transfer context are still lacking.
}
{
We describe the practical application of preconditioned Krylov methods to linear transfer problems of polarized radiation, possibly in a matrix-free context. The main aim is to clarify the advantages and drawbacks of various Krylov accelerators with respect to stationary iterative methods and direct solution strategies.
}
{
After a brief introduction to the concept of Krylov methods,
we report the convergence rate and the run time of various
Krylov-accelerated techniques combined with different formal solvers when applied to a 1D benchmark transfer problem of polarized radiation. 
In particular, we analyze the GMRES, BICGSTAB, and CGS Krylov methods, preconditioned with Jacobi, (S)SOR, or an incomplete LU factorization.
Furthermore, specific numerical tests were performed to study the
robustness of the various methods as the problem size grew.
}
{
Krylov methods accelerate the convergence, reduce the run time, and improve the robustness (with respect to the problem size)
of standard stationary iterative methods.
Jacobi-preconditioned Krylov methods outperform 
SOR-preconditioned stationary iterations  in all respects.
In particular, the Jacobi-GMRES method offers the best overall performance
for the problem setting in use.
}
{
Krylov methods can be more challenging to implement than stationary iterative methods. However, an algebraic formulation of the radiative transfer problem allows one to apply
and study Krylov acceleration strategies with little effort.
Furthermore, many available numerical libraries implement matrix-free Krylov routines, enabling an almost effortless transition to Krylov methods.}


\maketitle

\section{Introduction}\label{sec:intro}

When looking for numerical solutions to transfer problems
of polarized radiation,
it is common to rely on stationary iterative methods.
In particular, fixed point iterations with Jacobi, block-Jacobi, or Gauss-Seidel preconditioning have been successfully employed.
An illustrative convergence analysis of the stationary iterative methods usually used in the numerical transfer of polarized radiation is presented 
in the first paper of this series \citep{janett2021a}.
Unfortunately, stationary iterative methods often show unsatisfactory convergence rates when applied to large problems. By contrast, Krylov iterative methods gained popularity in the last decades, proving to be highly effective solution strategies, especially when dealing with large and sparse linear systems.

Among the various Krylov techniques, the biconjugate gradient stabilized (BICGSTAB) and the generalized minimal residual (GMRES) methods
are very common choices, especially when dealing with nonsymmetric systems \citep{meurant2020}.
As for stationary iterative methods, suitable preconditioning can significantly improve the convergence of Krylov iterations.

In the radiative transfer context, both BICGSTAB and GMRES methods have been employed for various applications,
such as
time-dependent fluid flow problems \citep{klein1989,hubeny2007new},
problems arising from finite element discretizations \citep{castro2015spatial,badri2019preconditioned},
and nonlinear radiative transfer problems for cool stars \citep{lambert2015new}.
In particular, Krylov methods have already been employed to linear transfer problems of unpolarized radiation, showing promising results.
We would like to mention the first convergence studies on the BICGSTAB method
by \citet{paletou2009conjugate} and \citet{anusha2009preconditioned}.

\citet{nagendra2009} first applied BICGSTAB to the transfer of polarized radiation in 1D atmospheric models.
\citet{anusha2011b} and~\citet{anusha2011c} employed the same method 
in 2D and 3D geometries, while \citet{anusha2012e}
additionally included angle-dependent
partial frequency redistribution (PRD) effects.
For the sake of simplicity, these papers are limited
to the study of hypothetical lines in isothermal atmospheric models.
Finally, \citet{anusha2013} used the BICGSTAB method to model
the Ca~{\sc ii} K 3993\,{\AA} resonance line
using an ad-hoc 3D atmospheric model,
and \citet{sampoorna2019} applied the GMRES method
to model the D$_2$ lines of Li~{\sc i} and Na~{\sc i},
taking the hyperfine structure of these atoms into account.
However, Krylov methods have not been fully exploited in the numerical solution of transfer of polarized radiation.
Indeed, only a few applications of Krylov methods aim to model the polarization profiles in realistic atmospheric models.
Moreover, numerical studies investigating multiple preconditioned Krylov methods, in terms of convergence and run time, and presenting key implementation details, are still lacking in this context.

This article is organized as follows. Section~\ref{sec:krylov} recalls the idea that lies behind Krylov iterative techniques, with a special focus on the GMRES and BICG methods. Section~\ref{sec:CRD_problem} presents a benchmark analytical problem,
its discretization, and its algebraic formulation,
and shows how to apply Krylov methods in the radiative transfer context.
In Section~\ref{sec:experiments}, we describe the problem settings and
the solvers options and present a quantitative comparison between
Krylov and stationary iterative methods
with respect to convergence, robustness, and run time.
Finally, Section~\ref{sec:conclusions} provides remarks and conclusions,
which are also generalized to more complex problems.

\section{Krylov methods}\label{sec:krylov}
This section contains a gentle, but not exhaustive, introduction to Krylov methods.
For more interested readers, \citet{ipsen1998idea} explain the idea behind Krylov methods, while \citet{saad2003iterative} and \citet{meurant2020} 
provide a comprehensive and rigorous discussion on the topic.

We consider the linear system
\begin{equation}\label{linear_system}
 A\xx=\bb,
\end{equation}
where the nonsingular matrix $A\in\mathbb R^{N\times N}$ and the vector $\bb \in \mathbb R^N$ are given, and the solution
$\xx  \in \mathbb R^N$ is to be found.
Given an initial guess $\xx^0$, a Krylov method approximates
the solution $\xx$ of~\eqref{linear_system} with the sequence of vectors $\xx^n\in\mathcal{K}_n(A,\bb)$, with $n=1,...,N$, where
\begin{equation}\label{eq:krylov}
\mathcal{K}_n(A,\bb) = \text{span}\left\{\bb,A\bb, A^2\bb,\ldots,A^{n-1}\bb\right\}
\end{equation}
is the $n$th Krylov subspace
generated by $\bb$. This is the case for the (very common) initial guess $\xx^0 = \mathbf{0}$. In general, $\xx^n-\xx^0\in\mathcal{K}_n(A,\rr^0)$, where $\rr^0 = \bb - A\xx^0$ is the initial residual.
 
 We first try to clarify why it is convenient to look for an approximate solution $\xx^n$ of the linear system~\eqref{linear_system} inside the Krylov subspace $\mathcal{K}_n(A,\bb)$.
The minimal polynomial of $A$ is the monic polynomial $q$ of minimal degree such that $q(A)=0$. If $A$ has $d$ distinct eigenvalues $\lambda_1,\ldots,\lambda_d$, its $m$th degree minimal polynomial reads
\begin{equation}\label{eq:min_pol}
  q(t) = \prod_{j=1}^d (t - \lambda_j)^{m_j}, 
\end{equation}
where $m_j$ is the multiplicity\footnote{The index $m_j$ is equal to the size of the largest Jordan block associated with $\lambda_j$. It is possible to show that $m_j \leq m_j^{A} - m_j^{G} + 1$, where $m_j^{A}$  and $m_j^{G}$  are the algebraic and geometric multiplicities associated with $\lambda_j$. If $A$ is diagonalizable,  for all $j$ we have $m_j^{A}=m_j^{G}$ and therefore $m_j=1$ for all $j$.} of  $\lambda_j$ (that is, the $j$th root of $q$) and $m = \sum_{j=1}^d m_j$. 
Polynomial~\eqref{eq:min_pol} can be rewritten in the form
\begin{equation*}
    q(t) = \sum_{j=0}^m\alpha_jt^j, 
\end{equation*}
where $\alpha_1,\ldots,\alpha_m$ are unspecified scalars and $\alpha_0 = \prod_{j=1}^d(-\lambda_j)^{m_j}$.
Evaluating the polynomial in $A$, one obtains
\begin{equation*}
    q(A) = \alpha_0 I\hspace{-0.1em}d + \alpha_1A + \ldots + \alpha_m A^m = 0,
\end{equation*}
which allows us to 
express the inverse $A^{-1}$ in terms of powers of $A$, namely,
\begin{equation}\label{eq:A_inv}
    A^{-1} = -\frac{1}{\alpha_0}\sum_{j=0}^{m-1}\alpha_{j+1} A^j,
\end{equation}
showing that the solution vector $\xx=A^{-1}\bb$ belongs to the Krylov subspace $\mathcal{K}_m(A,\bb)$. Hence, a smaller degree
of the minimal polynomial of $A$ (e.g., $m\ll N$)
generally corresponds to a faster convergence of Krylov methods, since $\xx\in\mathcal{K}_m(A,\bb)\subset \mathbb{R}^N$.
We remark that~\eqref{eq:A_inv} is well defined if $A$ is nonsingular:
 if $\lambda_j\neq 0$ for all $j$ then $\alpha_0 \neq 0$.

Secondly, we try to better characterize 
the smallest Krylov space containing $\xx$.
The minimal polynomial of $\bb$ with respect to $A$ is the monic polynomial $p$ of minimal degree such that $p(A)\bb=0$.
Denoting with $g$ the degree of $p$,
it is possible to show that $\mathcal{K}_g(A,\bb)$
is invariant under $A$, that is
\begin{equation*}
    \text{dim}\left(\mathcal{K}_n(A,\bb)\right)=\min\{n,g\},
\end{equation*}
meaning that $\mathcal{K}_n(A,\bb)=\mathcal{K}_g(A,\bb)$ for all $n\geq g$.
Therefore, a Krylov method applied to~\eqref{linear_system} will typically terminate in $g$ iterations, where $g$ can be much smaller than $N$. 
Since $g\leq m\leq N$, a Krylov method terminates after at most $N$ steps\footnote{The Cayley-Hamilton theorem implies $g \leq N$.}.
In practice, approximate solutions of system~\eqref{linear_system} are sufficient and
a suitable termination condition allows a Krylov method to converge in less than $g$ iterations.

Crucially, because of the structure of $\mathcal{K}_n(A,\bb)$, the computation of $\xx^n$ mainly requires a series of matrix-vector products involving $A$. This paradigm is particularly beneficial when $n$ is large and $A$ is sparse, in which cases matrix-vector products can be computed efficiently.
Moreover, Krylov methods are also handy when there is no direct access to the entries of the matrix $A$
and its action is only encoded in a routine that returns $A\vv$ from an arbitrary input vector $\vv$.
The various Krylov methods differ from one other in two main aspects.
The first difference is in the way they construct the subspace $\mathcal{K}_n(A,\bb)$.
In practice, the direct use of definition~\eqref{eq:krylov} is not convenient to build $\mathcal{K}_n(A,\bb)$, because the spanning vectors $\bb,\ldots,A^{n-1}\bb$ can become closer and closer being linearly dependent as $n$ grows, resulting in numerical instabilities.
The second difference is the way they select the ``best'' iterate $\xx^n\in\mathcal{K}_n(A,\bb)$.

The following sections introduce three common Krylov methods usually
applied to nonsymmetric linear systems: the GMRES, BICGSTAB and  conjugate gradient squared (CGS) methods.

\subsection{GMRES}\label{sec:gmres}
\cite{saad1986gmres} introduced one of the best known and mostly applied Krylov solvers:
the GMRES method. At iteration $n$, the GMRES method first constructs an orthonormal basis of $\mathcal{K}_n(A,\bb)$ by the Arnoldi algorithm, which is a modified version of the Gram-Schmidt procedure adapted to Krylov spaces.
Secondly, the GMRES method sets the approximate solution $\xx^n\in\mathcal{K}_n(A,\bb)$ minimizing the Euclidean norm of the residual $\rr^n = \bb - A\xx^n$, that is,
solving the least squares problem of size $n$ given by
\begin{equation*}
    \|\rr^n\|_2 = \| \bb - A\xx^n\|_2 =\min_{\zz\in\mathcal{K}_n}\| \bb - A\zz\|_2.
\end{equation*}
Assuming exact arithmetic, the GMRES method returns the exact solution in at most $g$ iterations. 
Regarding computational complexity,
the most expensive stages in the $n$th GMRES iteration is a $O(N^2)$
matrix-vector product involving $A$, possibly
decreasing to $O(N)$ for sparse $A$,
and a nested loop requiring $O(nN)$ floating point operations. 
In terms of memory, $n$ basis vectors have to be stored.
Therefore, the GMRES method becomes inefficient for large $n$.
A standard remedy to this problem is to employ restarts, that is, 
to terminate the GMRES method after $k$ iterations,
using $\xx^k$ as the initial guess for a subsequent GMRES run, until a desired termination condition is met. 
A suitable choice of $k$ can significantly reduce the time to solution, but convergence is not always guaranteed.
However, if $A$ is positive definite, GMRES converges for any $k\geq1$.

\subsection{CGS and BICGSTAB}\label{sec:cgs}
If $A$ is symmetric, the Arnoldi iteration can be simplified,
avoiding nested loops, obtaining the so-called Lanczos procedure, that requires just three vectors to be saved at once. 
%
For the nonsymmetric case, this procedure can be extended to the Lanczos biorthogonalization algorithm, that keeps in memory six vectors at once, yielding to
a significant
storage reduction
with respect to the Arnoldi method. On the other hand, two matrix-vector products (involving $A$ and $A^T$) are needed at each iteration.
Starting from two arbitrary vectors $\vv_1$ and $\ww_1$
satisfying $(\vv_1,\ww_1)=1$,
this algorithm builds two bi-orthogonal
bases $\{\vv_i\}_{i=1}^n$ and $\{\ww_i\}_{i=1}^n$ for $\mathcal{K}_n(A,\vv_1)$ and $\mathcal{K}_n(A^T,\ww_1)$, respectively.

In the $n$th iteration, the biconjugate gradient (BCG) method imposes
$\rr^n\perp\mathcal{K}_n(A^T,\ww_1)$ and $\xx_n \in\mathcal{K}_n(A,\vv_1)$. The BCG method efficiently solves two linear systems with coefficient matrices $A$ and $A^T$,
whose residual expressions are given by
\begin{equation}\label{eq:bcg_expressions}
\begin{aligned}
& \rr^n=p_n(A)\bb,\\
& \tilde{\rr}^n=p_n(A^T)\bb,
\end{aligned}
\end{equation}
$p_n$ being a $n$th-degree polynomial\footnote{From \eqref{eq:krylov}, by definition, we have that $\xx^n=p_{n-1}(A)\bb$. Then $$\rr^n=\bb-A\xx^n=\bb-Ap_{n-1}(A)\bb =p_{n}(A)\bb$$.} and $\tilde{\rr}^n = \bb - A^T\xx_n$.

Since the solution of the linear system involving $A^T$ is generally not needed and $A^T$ could not even be explicitly available,
transpose-free variants are used.
Two well-known transpose-free variants are the CGS \citep{sonneveld1989cgs} and the BICGSTAB \citep{van1992bi} methods. These variants avoid the use of $A^T$ replacing the residual expressions~\eqref{eq:bcg_expressions}
by $\rr_{\text{CGS}}^n = p_n^2(A)\bb$ and $\rr^n_{\text{BICGSTAB}} = p_n(A)\tilde{p}_n(A)\bb$,
respectively, where $\tilde{p}_n$ is a $n$th-degree polynomial recursively defined.
Regarding complexity, both CGS and BICGSTAB, in each iteration, require two matrix-vector products involving $A$, plus $O(N)$ additional floating point operations.
Both CGS and BICGSTAB are less robust than GMRES, since their convergence can stagnate for pathological cases \citep[see, e.g.,][]{gutknecht1993variants}.
Moreover, CGS rounding errors tend to be more pronounced with respect to BICGSTAB, possibly resulting in an unstable convergence.

\subsection{Preconditioning}\label{sec_prec}
Preconditioning the linear system~\eqref{linear_system} with the nonsingular matrix $P\in\mathbb R^{N\times N}$, that is, considering the equivalent system
\begin{equation}\label{eq:prec_sys}
    P^{-1}A\xx=P^{-1}\bb,
\end{equation}
can significantly increase robustness and convergence speed of
stationary iterative methods.
This is also particularly true for Krylov solvers. 
Crucially, the action of $P^{-1}$ on an input vector has to be computationally inexpensive, while improving the conditioning of $P^{-1}A$ with respect to the one of $A$.
In practice, Krylov methods must perform the product $\ww= P^{-1}A\vv$ for an arbitrary vector $\vv$. This operation is naively performed in two stages: $\tilde{\ww} = A\vv$ followed by $\ww = P^{-1}\tilde{\ww}$. 
Since $P$ usually depends on $A$, it is sometimes possible to obtain the action of $P^{-1}A$ in a more efficient way, reducing the number of floating point operations.
Since Krylov methods do not require direct access
to the entries of the preconditioner $P$,
the action of $P^{-1}$ (or directly of $P^{-1}A$)
can be also provided as a routine.

We would like to present the most common choices for $P$
that are numerically investigated in Section~\ref{sec:experiments}.
Given the decomposition $A=D+U+L$, with $D$, $U$, and $L$ being respectively the diagonal and the strictly upper and lower triangular parts of $A$, and $0<\omega<2$, we consider
\begin{enumerate}\setlength\itemsep{0.8em}
 \item $P=D$: Jacobi method;
 \item $P=\omega^{-1} D+L$ or $P=\omega^{-1}D+U$: successive over relaxation (SOR) method;
 \item $P=\omega(2-\omega)^{-1}\tilde{L}\tilde{U}$ with  $\tilde{L}=\omega^{-1}D+L$ and $ \tilde{U} = D^{-1}(\omega^{-1}D+U)$: symmetric SOR method (SSOR);
\item $P=\tilde{L}\tilde{U}$, where $\tilde{L}$ and $\tilde{U}$ are lower and upper triangular matrices approximating the factors of the LU factorization of $A$: incomplete LU factorization (ILU).
\end{enumerate}
The Jacobi method is numerically appealing since the action of $P^{-1}$
is cheaply computed (in parallel).
The SOR method usually provides faster convergence than Jacobi, at the price of $P^{-1}$ being more expensive to apply (only sequentially).
The SSOR method typically reduces the run time of SOR, especially if there is no preferential direction in the coupling between unknowns.
Moreover, a suitable tuning of $\omega$ significantly improves convergence, as shown in \citet{janett2021a} 
for stationary iterative methods. We remark that, for $\omega = 1$, SOR reduces to the Gauss-Seidel (GS) method and SSOR reduces to symmetric Gauss-Seidel.
ILU preconditioners are obtained trough an inexact Gaussian elimination, where some elements of the LU factorization are dropped. Typically, $\tilde{L}$ and $\tilde{U}$ are constructed such that their sum has the same sparsity pattern of $A$ (a.k.a. ILU with no fill-in or ILU(0)). As in the (S)SOR case, ILU preconditioners are sequential, but variants capable of extracting some parallelism can be applied.  

Finally, the Jacobi and (S)SOR methods can be encoded in matrix-free routines.
By contrast, to apply an ILU preconditioner, the matrices $A$, $\tilde{L}$ and $\tilde{U}$ must be assembled and stored.

\section{Benchmark problem}\label{sec:CRD_problem}
In this section, we present the continuous formulation of a 
benchmark linear transfer problem of polarized radiation. 
We then consider its discretization and algebraic formulation.
Finally, we remark on how the matrix-free Krylov methods presented in Sections~\ref{sec:gmres} and~\ref{sec:cgs} are applied in this context.
The reader is referred to \citet[][Sect.~3]{janett2021a} for a more detailed presentation of the problem,
along with its discretization and algebraic formulation. 

The problem is formulated within the framework of the complete frequency redistribution (CRD) 
approach presented in \citet{landi_deglinnocenti+landolfi2004}.
We consider a two-level atom with total angular momenta $J_u = 1$ and $J_\ell = 0$,
with the subscripts $\ell$ and $u$ indicating the lower and upper level, respectively.
Since $J_\ell=0$, the lower level is 
unpolarized by definition.
Stimulated emission it is not considered here, since it is completely negligible in solar applications. 
For the sake of simplicity, 
continuum processes are neglected, as are magnetic and bulk velocity fields.
A homogeneous and isothermal, one-dimensional (1D) plane-parallel atmosphere is considered.
Under the aforementioned assumptions, the spatial dependency of the physical quantities entering the problem is fully described by the height coordinate $z$. 
Moreover, the problem is characterized by cylindrical symmetry around the vertical
and, consequently, the angular dependency of the problem variables is fully described by the inclination $\theta$ with
respect to the vertical, or equivalently by $\mu=\cos(\theta)$.
\subsection{Continuous problem}
The physical quantities entering the problem 
are in general functions of the 
height $z\in[z_{\min},z_{\max}]$,
the frequency $\nu\in[\nu_{\min},\nu_{\max}]$, 
and the propagation direction of the radiation ray under consideration,
identified by $\mu\in[-1,1]$.

Due to the cylindrical symmetry of the problem, the only nonzero components of the radiation 
field tensor are $\bar{J}^0_0$ and $\bar{J}^2_0$. Consequently, the only nonzero multipolar components of the source 
function are $\sigma^0_0$ and $\sigma^2_0$.
Setting the angle $\gamma=0$ in the expression of the polarization tensor $\mathcal{T}^2_{0,i}$, one finds that
the only nonzero source functions 
are
\begin{align}
S_{\!I}(z,\mu) & = \sigma_0^0(z)+\mathcal{T}^2_{0,1}(\mu) \sigma_0^2(z), \label{sourceI}\\
S_{\!Q}(z,\mu) & = \mathcal{T}^2_{0,2}(\mu) \sigma_0^2(z), \label{sourceQ}
\end{align}
where $\mathcal{T}^2_{0,1}(\mu) = \sqrt{2}(3\mu^2-1)/4$ and $\mathcal{T}^2_{0,2}(\mu) = \sqrt{2}(3\mu^2-3)/4$. 
The only nonzero Stokes parameters are therefore $I$ and $Q$. Their propagation 
is described by the decoupled differential equations
\begin{align}\label{RTE_delo1}
	\mu\frac{\mathrm{d}}{\mathrm{d}z}
	I(z,\mu,\nu) &= - \eta(z,\nu)
	\left[I(z,\mu,\nu) - S_{\!I}(z,\mu) 
	\right],\\\label{RTE_delo2}
	\mu\frac{\mathrm{d}}{\mathrm{d}z}
	Q(z,\mu,\nu) &= - \eta(z,\nu)
	\left[Q(z,\mu,\nu) - S_{\!Q}(z,\mu) 
	\right],
\end{align}
where $\eta$ is the absorption coefficient. 
The initial 
conditions are given by
\begin{align*}
& I(z_{\min},\mu,\nu)=1,  &&\text{ if $\mu>0$},\\
& Q(z_{\min},\mu,\nu)=0,  &&\text{ if $\mu>0$},\\
& I(z_{\max},\mu,\nu) = Q(z_{\max},\mu,\nu)=0, &&\text{ if $\mu<0$}.
\end{align*}
In order to linearize the problem, we assume that the absorption coefficient 
$\eta$ is known a priori and fixed or, equivalently, we discretize the atmosphere through a 
fixed grid in frequency-integrated optical depth \citep[see][Sect.~3.2]{janett2021a}.

The explicit expressions of $\bar{J}^0_0$ and $\bar{J}^2_0$ are
\begin{align}
\bar J^0_0(z)  = & \int {\rm d} \nu\frac{\phi(\nu)}{2} 
	\oint \mathrm{d}\mu\,
	I(z,\mu,\nu), \label{J00}\\
	\bar J^2_0(z) = & \int {\rm d} \nu\frac{\phi(\nu)}{2}
	\oint \mathrm{d} \mu
	\left[\mathcal{T}^2_{0,1}(\mu)I(z,\mu,\nu) + \mathcal{T}^2_{0,2}(\mu)Q(z,\mu,\nu)\right],\label{J02}
\end{align}
where the damping constant entering the absorption profile $\phi$ is fixed to $a=10^{-3}$.
Finally, the explicit expressions of $\sigma^0_0$ and $\sigma^2_0$ are 
\citep[see][Eq.~(18)]{janett2021a}
\begin{align}
\sigma_0^0(z)& = \xi\bar J^0_0(z) + \epsilon\notag\\
	&= \xi\!\int {\rm d} \nu\frac{\phi(\nu)}{2}
	\oint \mathrm{d} \mu\,
	I(z,\mu,\nu) + \epsilon, \label{S00}\\
	\sigma_0^2(z) &= \xi\bar J^2_0(z)\nonumber\\
	&= \xi\!\int{\rm d}\nu\frac{\phi(\nu)}{2}
	\oint \mathrm{d} \mu
	\left[\mathcal{T}^2_{0,1}(\mu)I(z,\mu,\nu) + \mathcal{T}^2_{0,2}(\mu)Q(z,\mu,\nu)\right],\label{S02}
\end{align}
where we set the thermalization parameter to $\epsilon=10^{-4}$
and $\xi=1-\epsilon$.
In the equations above,
we assumed a negligible depolarizing rate of the upper level due 
to elastic collisions $\delta_u^{(K)}=0$, and set the Planck function
(in the Wien limit) to $W = 1$.
We finally recall that for the considered atomic model 
$w_{J_u J_\ell}^{(0)} = w_{J_u J_\ell}^{(2)} = 1$.
\begin{figure*}
    \centering
    \includegraphics[width=\textwidth]{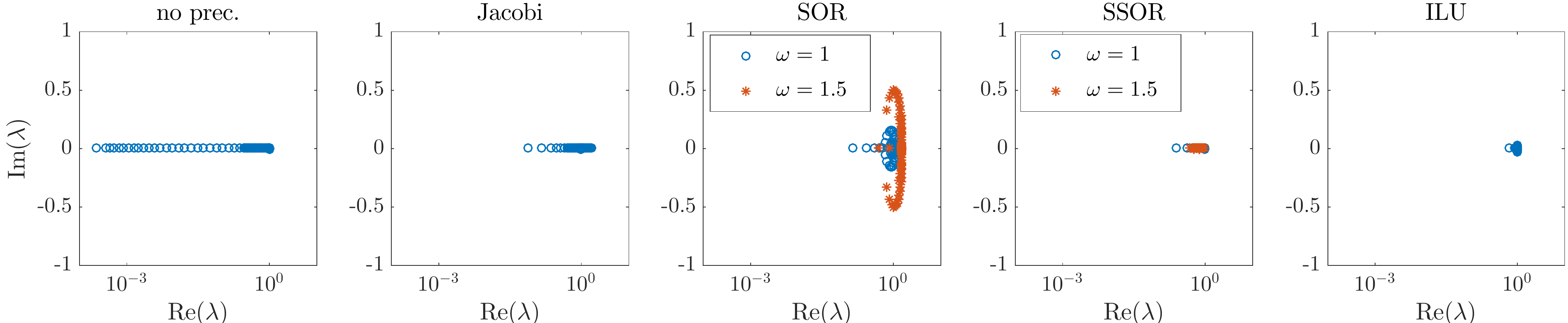}
    \caption{Spectrum of $P^{-1}A$ on the complex plane
    for different preconditioners
    for the DELO-linear formal solver, with $N_s = 80$ and $N_\mu=N_\nu=20$.  From here on
    the abbreviation ``no prec.'' indicates that no preconditioner is used, that is, $P=I\hspace{-0.1em}d$.}
    \label{fig:eigs}
\end{figure*}
\subsection{Discrete problem}
We discretize the continuous variables
$z$, $\nu$, and $\mu$ with the numerical grids
$\{z_k\}_{k=1}^{N_s}$, $\{\mu_m\}_{m=1}^{N_\mu}$, and $\{\nu_p\}_{p=1}^{N_{\nu}}$, respectively.
The quantities of the problem are now approximated at the nodes only.
The discrete version of equations~\eqref{sourceI}--\eqref{sourceQ}
can be expressed in matrix form
\begin{equation}\label{matricial_form_2}
    \mathbf{S}=T\pmb{\sigma},\;\;\text{ with }
T\in\mathbb R^{2 N_s N_\mu \times2 N_s},
\end{equation}
where $\mathbf{S}\in\mathbb R^{2 N_s N_\mu }$ and $\pmb{\sigma}\in\mathbb R^{2 N_s}$
collect the discretized source functions 
$S_{\!I}(z_k,\mu_m)$ and $S_{\!Q}(z_k,\mu_m)$,
and the multipolar components of the source function $\sigma_0^0(z_k)$ and $\sigma_0^2(z_k)$, respectively.
The entries of the matrix $T$ depend on the coefficients
appearing in \eqref{sourceI}--\eqref{sourceQ}. 

Similarly, the transfer equations~\eqref{RTE_delo1}--\eqref{RTE_delo2}
are expressed in matrix form
\begin{equation}\label{matricial_form_3}
    \mathbf{I}=\Lambda\mathbf{S}+\mathbf{t},\;\;\text{ with }
\Lambda\in\mathbb R^{2  N_s N_\mu N_\nu \times 2 N_s N_\mu },
\end{equation}
where $\mathbf{I}\in\mathbb R^{2 N_s  N_\mu N_\nu}$ collects
the discretized Stokes parameters $I(z_k,\mu_m,\nu_p)$ and $Q(z_k,\mu_m,\nu_p)$, while $\mathbf{t}\in\mathbb R^{2 N_s  N_\mu N_\nu}$ represents the radiation transmitted from the boundaries.
The entries of the matrix $\Lambda$ depend
on the numerical method (a.k.a. formal solver) used to solve the transfer equations~\eqref{RTE_delo1}--\eqref{RTE_delo2}, on the spatial grid, and on the eventual numerical conversion to the optical depth scale.

Finally,
the matrix version of~\eqref{S00}--\eqref{S02} reads
\begin{equation}\label{matricial_form_1}
\pmb{\sigma}=J\mathbf{I}+\mathbf{c},\;\;\text{ with }
J\in\mathbb R^{2 N_s\times2 N_s N_\mu N_\nu},
\end{equation}
where $\pmb{\sigma} = [\sigma_0^0(z_1),\sigma_0^2(z_1),\sigma_0^0(z_2),\ldots,\sigma_0^2(z_{N_s})]^T$ and, accordingly, $\mathbf{c} = [\epsilon, 0, \epsilon,0,\ldots,\epsilon, 0]^T$. The matrix $J$ depend on the choice of the angular and spectral quadratures
used in~\eqref{S00}--\eqref{S02}.
\citet{janett2021a}
provide an explicit description of the matrices $T$, $\Lambda$, and $J$ for a more general radiative transfer setting.

By choosing $\pmb{\sigma}$ as the unknown vector,
Equations~\eqref{matricial_form_2}--\eqref{matricial_form_1} can then be combined
in the linear system given by
\begin{equation}\label{crd_linear_system}
 (I\hspace{-0.1em}d-J\Lambda T)\pmb{\sigma}=J\mathbf{t}+\mathbf{c} \quad \iff \quad A\xx=\bb,
\end{equation}
with $I\hspace{-0.1em}d-J\Lambda T$ being a square matrix of size $2N_s$.
Referring to \eqref{linear_system}, we set $A = I\hspace{-0.1em}d-J\Lambda T$, $\xx=\pmb{\sigma}$,
and $\bb = J\mathbf{t}+\mathbf{c}$. 
Alternatively, it is also possible to consider $\mathbf{I}$ as unknown of the problem,
obtaining the following linear system
\begin{equation}\label{crd_linear_system2}
 (I\hspace{-0.1em}d-\Lambda TJ)\mathbf{I}=\Lambda T\mathbf{c}+\mathbf{t},
\end{equation}
with $I\hspace{-0.1em}d-\Lambda TJ$ being a block sparse square matrix of size $2N_sN_\mu N_\nu$
with dense blocks.

The linear problem~\eqref{crd_linear_system} is considered in the following sections,
because of its smaller size in the current setting.
However, the formulation~\eqref{crd_linear_system2} is likely to be more favorable
to Krylov methods, since the corresponding linear system is both larger and sparse. 
In both cases the coefficient matrix $A$ is nonsymmetric and positive definite (cf. Figure~\ref{fig:eigs}). 

\subsection{Matrix-free Krylov approach}
The explicit assembly of $A = I\hspace{-0.1em}d-J\Lambda T$ is convenient if the linear system~\eqref{crd_linear_system} must be solved multiple times, for example in the case of extensive numerical tests with different right-hand sides.
However, it is in general expensive to assemble $A$,
especially for large problems.
However, Krylov methods can also be applied in a matrix-free context,
where the action of $A$ is encoded in a routine and no direct access to its entries is required. For this reason, we provide in Algorithm~\ref{algo} a matrix-free routine to compute $\pmb{\sigma}_\text{out} =A\pmb{\sigma}_\text{in}$ for any arbitrary input vector $\pmb{\sigma}_\text{in}$.
Algorithm~\ref{algo} is modular, allowing preexisting radiative transfer routines to be reused. 

We remark that it is possible to assemble $A$ column-by-column,
obtaining its $j$th column applying Algorithm~\ref{algo}
to a point-like source function  $\pmb{\sigma}_i=\delta_{ij}$ for $i,j = 1,\ldots,2N_s$.
\begin{algorithm} \label{algo}
\SetAlgoLined
\KwIn{$\pmb{\sigma}_\text{in}$, that is, $\{ \sigma_0^0(z_k)\}_{k=1}^{N_s}$ and $\{ \sigma_0^2(z_k)\}_{k=1}^{N_s}$}
\tcp{compute $\mathbf{S} = T\pmb{\sigma}_\text{in}$}
\For {$k=1,\ldots,N_s$}{
 \For {$m=1,\ldots,N_\mu$}{
 compute $S_{\!I}(z_k,\mu_m)$ and $S_{\!Q}(z_k,\mu_m)$ from \eqref{sourceI}--\eqref{sourceQ};
 }
 }
 \tcp{compute $\mathbf{I}=\Lambda\mathbf{S}$ with a formal solver}
  \For {$m=1,\ldots,N_\mu$}{
   \For {$p=1,\ldots,N_\nu$}{
    given $S_{\!I}(z_k,\mu_m)$ and $S_{\!Q}(z_k,\mu_m)$, compute $\{I(z_k,\mu_k,\nu_p)\}_{k=1}^{N_s}$ and $\{Q(z_k,\mu_k,\nu_p)\}_{k=1}^{N_s}$ by solving transfer equations \eqref{RTE_delo1} and \eqref{RTE_delo2} with a predefined formal solver and zero initial conditions (i.e., with $I(z_{\min},\mu,\nu) = 0$).
  }
 }
 \tcp{compute $\pmb{\sigma} = J\mathbf{I}$}
 \For {$k=1,\ldots,N_s$}{
given $\{I(z_k,\mu_k,\nu_p)\}_{k=1}^{N_s}$ and $\{Q(z_k,\mu_k,\nu_p)\}_{k=1}^{N_s}$,
compute $\bar J^0_0(z_k)$ and $\bar J^2_0(z_k)$ according to \eqref{J00}--\eqref{J02} using a predefined numerical quadrature.
 }
 \tcp{$\pmb{\sigma}_\text{out} = \pmb{\sigma}_\text{in} - \pmb{\sigma}$}
 \KwOut{$\pmb{\sigma}_\text{out} = \pmb{\sigma}_\text{in} - \xi\,[\bar J^0_0(z_1),\bar J^2_0(z_1),\ldots,\bar J^2_0(z_{N_s})]$}
 \caption{Compute $\pmb{\sigma}_\text{out} = A\pmb{\sigma}_\text{in}$, given $$\pmb{\sigma}_\text{in} = [\sigma_0^0(z_1),\sigma_0^2(z_1),\sigma_0^0(z_2),\ldots,\sigma_0^2(z_{N_s})].$$}

\end{algorithm}
\section{Numerical experiments} \label{sec:experiments}
In Section~\ref{sec:krylov}, we discussed the convergence properties,
complexity, and memory requirements of various Krylov methods.
In practice, the specific choice among the different methods is usually made according to run time, which is empirically tested.

In this section, different Krylov and stationary iterative methods
are applied to the benchmark linear problem \eqref{crd_linear_system}. In particular, we compare GMRES, CGS, and BICGSTAB Krylov methods, possibly preconditioned, and the Richardson, Jacobi, SOR, and SSOR stationary iterative methods.
The stationary iterative methods are described as preconditioned Richardson methods
of the form
%
\begin{equation}\label{eq:rich_iteration}
    \xx^{n+1} = \xx^n + P^{-1}(\bb - A\xx^n).
\end{equation}
\subsection{Discretization, quadrature, and formal solver}
We replaced the spatial variable $z$ by the 
frequency-integrated optical depth scale along the vertical $\tau$.
Moreover, we fixed the logarithmically spaced grid
\begin{equation*}
10^{-5} = \tau_1 < \tau_2 < \cdots < \tau_{N_s} = 10^4,
\end{equation*}
thus linearizing the problem.
%
We used $N_\mu$ Gauss-Legendre nodes and weights to discretize $\mu\in[-1,1]$ and compute the corresponding angular quadrature. 
The spectral line was sampled with $N_\nu$ frequency nodes
equally spaced in the reduced frequency interval $[-5,5]$ and
the trapezoidal rule was used as spectral quadrature.

A limited number of formal solvers was considered 
for the numerical solution of transfer equations~\eqref{RTE_delo1}--\eqref{RTE_delo2}:
the implicit Euler, DELO-linear, DELOPAR, and DELO-parabolic methods.
Further details on these
formal solvers are given by \citet{janett2017a,janett2017b,janett2018b}.
%
%
\begin{figure}
    \centering
    \includegraphics[width=0.46\textwidth]{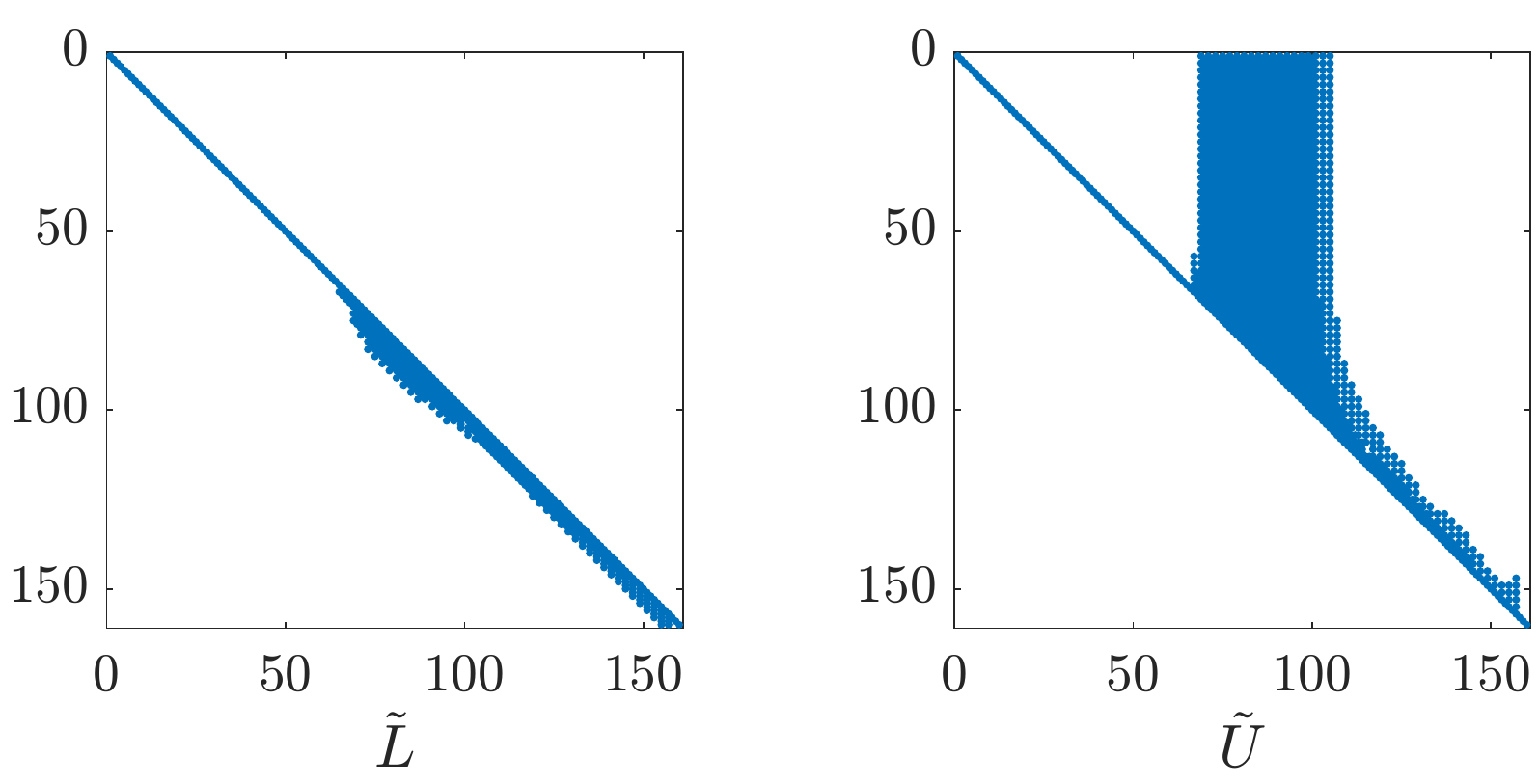}
    \caption{Sparsity pattern of $\tilde{L}$ and $\tilde{U}$
    from the ILU factorization of $A$
    for the DELO-parabolic formal solver,
    with $N_s=80$ and $N_{\mu}=N_{\nu}=20$.}
    \label{fig:sparsity_ILU}
\end{figure}
\subsection{Iterative solvers settings and implementation}
As initial guess for all the solvers we use $\sigma_0^0(z_k) = 1$ and $\sigma_0^2(z_k) = 0$ for all $k$, corresponding to $\pmb{\sigma}_{\text{in}}=[1,0,1, 0,\ldots,1,0]$.
Since multiple preconditioners are compared,
we use the following termination condition
\begin{equation*}
\|\rr^n\|_2/\|\bb\|_2<10^{-6},    
\end{equation*}
which is based on the unpreconditioned relative residual.

Between the two possible formulations of SOR presented in Section~\ref{sec_prec},
the structure of $A$ (cf. Figure~\ref{fig:sparsity_ILU}) suggests the use of $P=D + \omega U$,
which in fact provides a better convergence.
Moreover, the near-optimal damping parameter $\omega=1.5$ \citep[see][Sect. 2.3]{janett2021a}
is used for both SOR and SSOR preconditioners when applied to the Richardson method.
By contrast, no damping (that is, $\omega=1$) is used when
the SOR and SSOR preconditioners are applied to Krylov methods,
since no benefits are observed for different choices of $\omega$. 
The GMRES method is employed without restarts,
while the ILU factorization is applied with a threshold of $10^{-2}$
(a.k.a. ILUT\footnote{Replacing with zero nondiagonal elements of $\tilde{L}$ and $\tilde{U}$ if
$$ |\tilde{U}_{ij}| < 10^{-2} \|A_{* j}\|_2, \qquad  |\tilde{L}_{ij}| < 10^{-2} \|A_{* j}\|_2/\tilde{U}(j,j),$$ 
for $i,j=1,..,N$ and $i\neq j$.}).
Figure~\ref{fig:sparsity_ILU} exposes
an example of sparsity pattern of $\tilde{L}$ and $\tilde{U}$,
highlighting the most significant entries of $A$.

Krylov methods are applied to problem \eqref{crd_linear_system}
using MATLAB and its built-in functions: \texttt{gmres, bicgstab, cgs}
and \texttt{ilu}.
In practice, we explicitly assemble the matrix $A$
and store $P$ in sparse format,
using \texttt{mldivide} to apply $P^{-1}$ in~\eqref{eq:rich_iteration}.
The SSOR and ILU preconditioners
can be expressed as the product of lower and upper triangular matrices, that is,
in the form $P=\tilde{L}\tilde{U}$.
Since $P^{-1} = \tilde{U}^{-1}\tilde{L}^{-1}$,
the action of $P^{-1}$ can be conveniently computed
using sparse storage for $\tilde{L}$ and $\tilde{U}$.
We would like to mention that the Eisenstat's trick (not used here) can be conveniently applied to the SSOR preconditioning, reducing the run time by a factor up to two.

Matrix-free approaches are also adopted, 
using Algorithm~\ref{algo} and avoiding the explicit assembling of the matrix $A$. In this case, for performance reasons, we precompute the action of $P^{-1}$. Notice that the Krylov routines available in the most common numerical packages (e.g., MATLAB or PETSc) support a matrix-free format of $A$ and $P^{-1}$.
\begin{table}
\centering
    \begin{tabular}{l|c|c|c|c|c|c|c}
         $N_{\mu} = N_{\nu}$ & 20 & 30 & 40 & 50 & 60 & 70 & 80\\ 
         \hline 
         \rule{0pt}{1.1\normalbaselineskip}GMRES & 48 & 48 & 49 & 49 & 49 & 49 & 49 \\
         BICGSTAB & 58 & 57 & 55 & 58 & 60 & 57 & 57\\
         CGS & 54 & 55 & 54 & 57 & 55 & 55 & 55
    \end{tabular}
    \caption{Iterations to convergence for the GMRES, CGS, and BICGSTAB Krylov methods, with no preconditioning, for the DELO-linear formal solver, with $N_s = 40$,
    and varying $N_{\mu} = N_{\nu}$.}
    \label{tab:mu_nu_dependency}
\end{table}

\subsection{Convergence}
Figure~\ref{fig:eigs} displays the spectrum of $P^{-1}A$ on the complex plane
for different preconditioners. The clustering of the eigenvalues of $P^{-1}A$ corresponds to a better conditioning and, consequently, to a faster convergence of iterative methods when applied to~\eqref{eq:prec_sys}. Clustered eigenvalues also yield a reduction of the dimension of the Krylov space $\mathcal{K}_n(P^{-1}A,P^{-1}\bb)$. 
Table~\ref{tab:mu_nu_dependency} presents the number of iterations required
from different methods to converge
as a function of the number of nodes in the spectral and angular grids $N_{\mu}$ and $N_{\nu}$, respectively.
Analogously to \citet{janett2021a}, 
we observe that, once a minimal resolution is guaranteed, varying $N_{\mu}$ and $N_{\nu}$
has a negligible effect on the convergence behavior of the different Krylov solvers.
Therefore, the values $N_{\mu}=20$ and $N_{\nu}=20$ remain fixed in the following numerical experiments on convergence.

\begin{figure*}
    \centering
    \includegraphics[width=\textwidth]{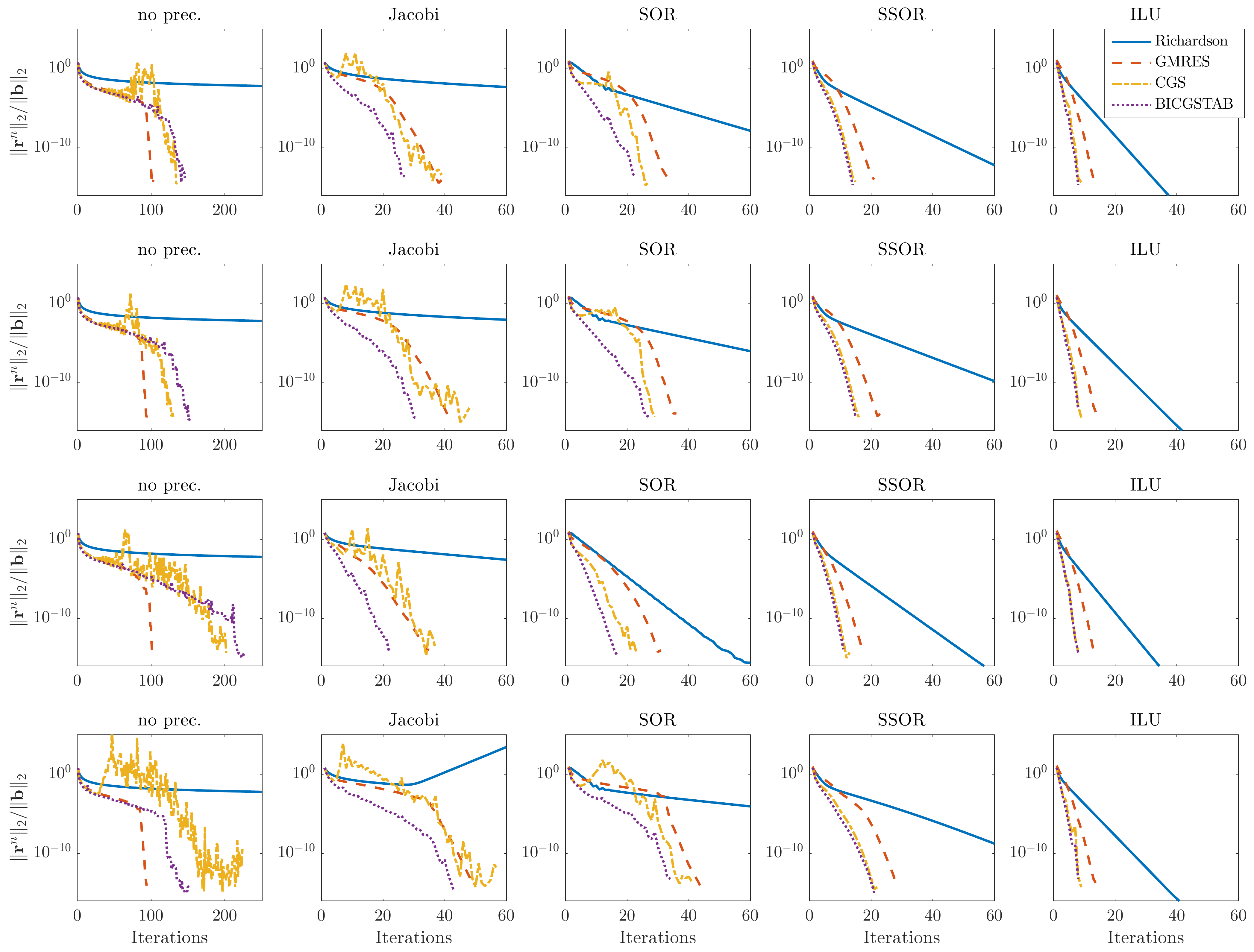}
    \caption{Relative residual, up to machine precision, vs iterations for various preconditioners and formal solvers with $N_s = 80$ and $N_\nu=N_\mu =20$. Each row corresponds to a different formal solver, from above: implicit Euler, DELO-linear, DELOPAR, DELO-parabolic. We note that a different scale is used in the horizontal axes of the first column, where no preconditioner is employed.}
    \label{fig:conv}
\end{figure*}

Figure~\ref{fig:conv} presents the convergence history of the Richardson and Krylov methods combined with different preconditioners and the formal solvers under investigation.
Krylov methods generally outperform stationary methods in terms of number of iterations to reach convergence. Among the Krylov methods,
preconditioned BICGSTAB and CGS methods result in the fastest convergence.
As predicted by Figure~\ref{fig:eigs}, SSOR and ILU preconditioners are the most effective.
On the other hand,
the use of a Jacobi preconditioner also seems ideal,
because of its overall fast convergence
and its cheap and possibly parallel application.

\begin{figure*}
        \centering
    \includegraphics[width=\textwidth]{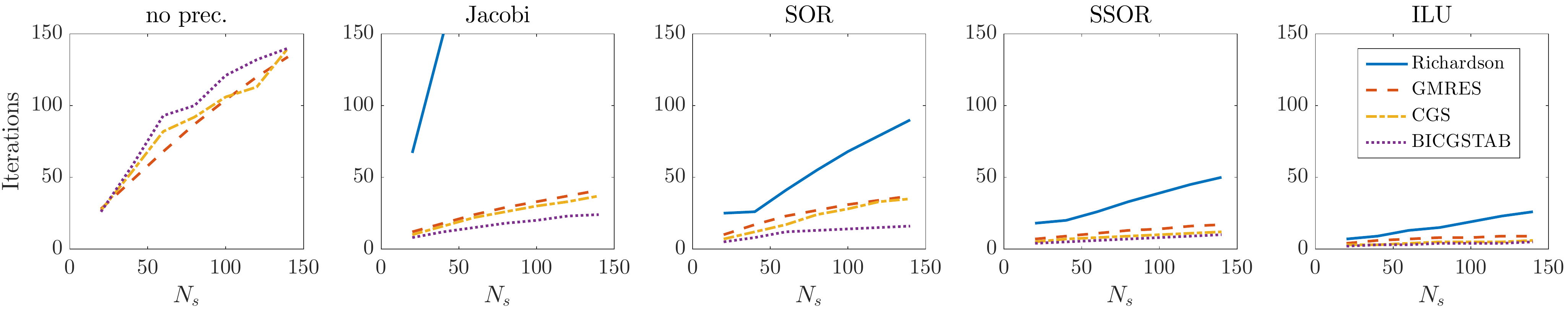}
    \caption{Graphical representation of Tables~\ref{tab:no_prec}--\ref{tab:ilu}.}
    \label{fig:conv2}
\end{figure*}

Tables~\ref{tab:no_prec}--\ref{tab:ilu} present
the number of iterations required by different methods to converge
as a function of $N_{s}$ for different preconditioners,
using the DELO-linear formal solver.
The same numerical experiments using the implicit Euler, DELOPAR,
and DELO-parabolic formal solvers show similar trends,
which are thus not shown.  With reference to Table~\ref{tab:no_prec}, we observe that the unpreconditioned Richardson method never converges in less than $10^4$ iterations.
Figure~\ref{fig:conv2} graphically represents
the content of Tables~\ref{tab:no_prec}--\ref{tab:ilu},
suggesting a superior scaling of Krylov methods with the problem size in terms of convergence. 
Since realistic radiative transfer problems are very large,
these results are particularly relevant for practical purposes.
\begin{table}[ht!]
    \centering
    \begin{tabular}{l|r|r|r|r|r|r|r}
         $N_s$ & \;\;20 & \;\;40 & \;\;60 & \;\;80 & 100 & 120 & 140\\ 
         \hline 
         \rule{0pt}{1.1\normalbaselineskip}Richardson & - & - & - & - & - & - & - \\
         GMRES & 28 & 48 & 68 & 87 & 104 & 120 & 134 \\
         BICGSTAB & 26 & 58 & 93 & 100 & 121 & 132  & 140 \\
         CGS      & 27 & 54 & 82 & 92 & 106 & 113 & 140
    \end{tabular}
    \caption{Iterations to convergence for the Richardson, GMRES, CGS, and BICGSTAB methods with no preconditioner for
    the DELO-linear formal solver,
    with $N_{\mu} =N_{\nu}=20$, and varying  $N_s$.}
    \label{tab:no_prec}
    \centering
    \begin{tabular}{l|r|r|r|r|r|r|r}
         $N_s$ & \;\;20 & \;\;40 & \;\;60 & \;\;80 & 100 & 120 & 140\\ 
         \hline 
         \rule{0pt}{1.1\normalbaselineskip}Richardson & 67 & 150 & 230 & 304 & 374 & 441 & 504 \\
         GMRES    & 12 & 18 & 24 & 29 & 33 & 37 & 41 \\
         BICGSTAB & 8  & 12 & 15 & 18 & 20 & 23 & 24 \\
         CGS      & 10 & 16 & 22 & 26 & 30 & 33 & 37
    \end{tabular}
    \caption{Same as Table~\ref{tab:no_prec},
    but for the Jacobi preconditioner.}
    \label{tab:jacobi}
    \centering
    \begin{tabular}{l|r|r|r|r|r|r|r}
         $N_s$ & \;\;20 & \;\;40 & \;\;60 & \;\;80 & 100 & 120 & 140\\ 
         \hline 
         \rule{0pt}{1.1\normalbaselineskip}Richardson & 25 & 26 & 41 & 55 & 68 & 79 & 90  \\
         GMRES    & 10 & 17 & 23 & 27 & 31 & 34 & 37  \\
         BICGSTAB & 5  & 8 & 12 & 13  & 14 & 15 & 16  \\
         CGS      & 7  & 12 & 17 & 24 & 28 & 33 & 35
    \end{tabular}
    \caption{Same as Table~\ref{tab:no_prec},
    but for the SOR preconditioner.}
    \label{tab:sor}
    \centering
    \begin{tabular}{l|r|r|r|r|r|r|r}
         $N_s$ & \;\;20 & \;\;40 & \;\;60 & \;\;80 & 100 & 120 & 140\\ 
         \hline 
         \rule{0pt}{1.1\normalbaselineskip}Richardson & 18 & 20 & 26 & 33 & 39 & 45 & 50  \\
         GMRES    & 7 & 9 & 11 & 13 & 14 & 16 & 17  \\
         BICGSTAB & 4 & 5 & 6 & 7  & 8 & 9 & 10  \\
         CGS      & 5 & 7 & 8 & 9 & 10 & 11 & 12
    \end{tabular}
    \caption{Same as Table~\ref{tab:no_prec},
    but for the SSOR preconditioner.}
    \label{tab:ssor}
    \centering
    \begin{tabular}{l|r|r|r|r|r|r|r}
         $N_s$ & \;\;20 & \;\;40 & \;\;60 & \;\;80 & 100 & 120 & 140\\ 
         \hline 
         \rule{0pt}{1.1\normalbaselineskip}Richardson & 7 & 9 & 13 & 15 & 19 & 23 & 26 \\
         GMRES    & 4 & 6 & 7 & 8 & 8 & 9 & 9 \\
         BICGSTAB & 2 & 3 & 3 & 4 & 4 & 4 & 5 \\
         CGS      & 3 & 3 & 4 & 5 & 5 & 5 & 6 
    \end{tabular}
    \caption{Same as Table~\ref{tab:no_prec},
    but for the ILU preconditioner.}
    \label{tab:ilu}
\end{table}

\subsection{Run times}

Since the use of different formal solvers would not change the essence of the results,
this section only considers the use of the DELO-linear formal solver.

We first consider the case of $A$ being assembled explicitly.
Table~\ref{tab:run_times_assemble} reports the time spent to assemble $A$ using Algorithm~\ref{algo} repeatedly, varying the discretization parameters. We remark that Algorithm~\ref{algo} can be optimized for this context, where a point-like input is used. In fact, the computational complexity of the formal solution of \eqref{RTE_delo1}--\eqref{RTE_delo2} with a zero initial condition and a point-like source term can be largely reduced. For example, given $j\in\{1,\ldots,N_s\}$ and a point-like $\pmb{\sigma}$ in $z_j$, that is, $\pmb{\sigma}_i=\delta_{ij}$, the solution of \eqref{RTE_delo1}--\eqref{RTE_delo2} for $\mu<0$ (resp. $\mu>0$) is identically zero for all $z<z_j$ (resp. $z>z_j$) and it is given by an exponential decay for $z>z_j$ (resp. $z<z_j$).
Once $A$ is explicitly assembled, problem~\eqref{crd_linear_system}
can be solved with a negligible run time using, for example, a LU direct solver (taking 32 milliseconds for $N_s=500$, cf. Table~\ref{tab:run_times_assemble}).
We remark that the solve time is almost negligible, with respect to the time spent on assembly,
for all the iterative methods under consideration (besides unpreconditioned Richardson).
Therefore, Tables~\ref{tab:run_times_mf}-\ref{tab:run_times_mf2}-\ref{tab:run_times_mf3} report run times corresponding to a matrix-free approach, for which no time is spent to assemble $A$. With reference to Table~\ref{tab:run_times_mf}, the time spent to compute $P$ is of approximately 0.05, 4 and 9 seconds for the Jacobi, SOR, and SSOR preconditioners, respectively. The run time of the Richardson method with no preconditioning is not reported, because it exceeds one hour. Similarly, concerning Table~\ref{tab:run_times_mf2}, the time spent to compute $P$ is of approximately 0.2, 55 and 111 seconds for the Jacobi, SOR, and SSOR preconditioners, respectively.
In the matrix-free context, GMRES turns out to be the fastest Krylov method, since it requires only one application of $A$ for each iteration (as discussed in Section~\ref{sec:krylov}). Consider that the application of $A$, obtained using Algorithm~\ref{algo}, is computationally expensive.

Figure~\ref{fig:bar_plot} displays the total run times of the most relevant approaches,
highlighting the time spent to assemble $A$, to precompute the action of $P^{-1}$,
and to solve the linear system.
Accordingly, the most convenient methods are, assembling $A$ and then apply an LU direct solver or employing a Jacobi-GMRES matrix-free solver.

The best choice between the two approaches mostly depends on the problem size,
memory limitations, and accuracy requirements.
To this regard, Table~\ref{tab:run_times_mf3} reports how the run time of Jacobi-GMRES scales
with respect to the problem size (in terms of $N_s,N_\mu$ and $N_\nu$).
The comparison between Table~\ref{tab:run_times_assemble} and Table~\ref{tab:run_times_mf3} suggests that the matrix-free approach scales better as $N_s$ grows. On the other hand, it becomes advantageous to assemble $A$ as $N_\mu$ and $N_\nu$ grow (cf. Figure~\ref{fig:bar_plot2}). To explain this difference, we would like to remark that, in general, the nested loop over $N_\mu$ and $N_\nu$, containing a formal solver call over $N_s$, is the most expensive portion of Algorithm~\ref{algo}.
On the other hand, when Algorithm~\ref{algo} is used to assemble $A$, the computational complexity of the formal solution of \eqref{RTE_delo1}--\eqref{RTE_delo2} can be largely reduced. In this case, the computational bottleneck of Algorithm~\ref{algo}
is the calculation of $\bar J^0_0(z_k)$ and $\bar J^2_0(z_k)$, for $k=1,\ldots,N_s$.

We finally remark that both the selected methods can be applied in parallel,
since each loop in Algorithm~\ref{algo} can be replace by a ``parallel for''
and the application of the Jacobi preconditioner is embarrassingly parallel. 
\section{Conclusions}\label{sec:conclusions}

We described how to apply state-of-the-art Krylov methods to linear radiative transfer problems of polarized radiation. In this work, we considered the 
same benchmark problem as in \citet{janett2021a}.
However, Krylov methods can be applied to a wider range of settings,
such as the unpolarized case, two-term atomic models 
including continuum contributions, 3D atmospheric models with arbitrary magnetic and bulk velocity fields, and theoretical frameworks accounting for PRD 
effects.

In particular, we presented the convergence behavior and the run time measures of the GMRES, BICGSTAB, and CGS Krylov methods. Krylov methods 
significantly accelerate standard stationary iterative methods in terms of convergence rate, time-to-solution, and are favorable in terms of scaling 
with respect to the problem size.
Contrary to stationary iterative methods, Krylov-accelerated routines always converge (even with no preconditioning) for the considered numerical 
experiments.

\begin{figure*}[ht!]
    \centering
    \includegraphics[width=0.95\textwidth]{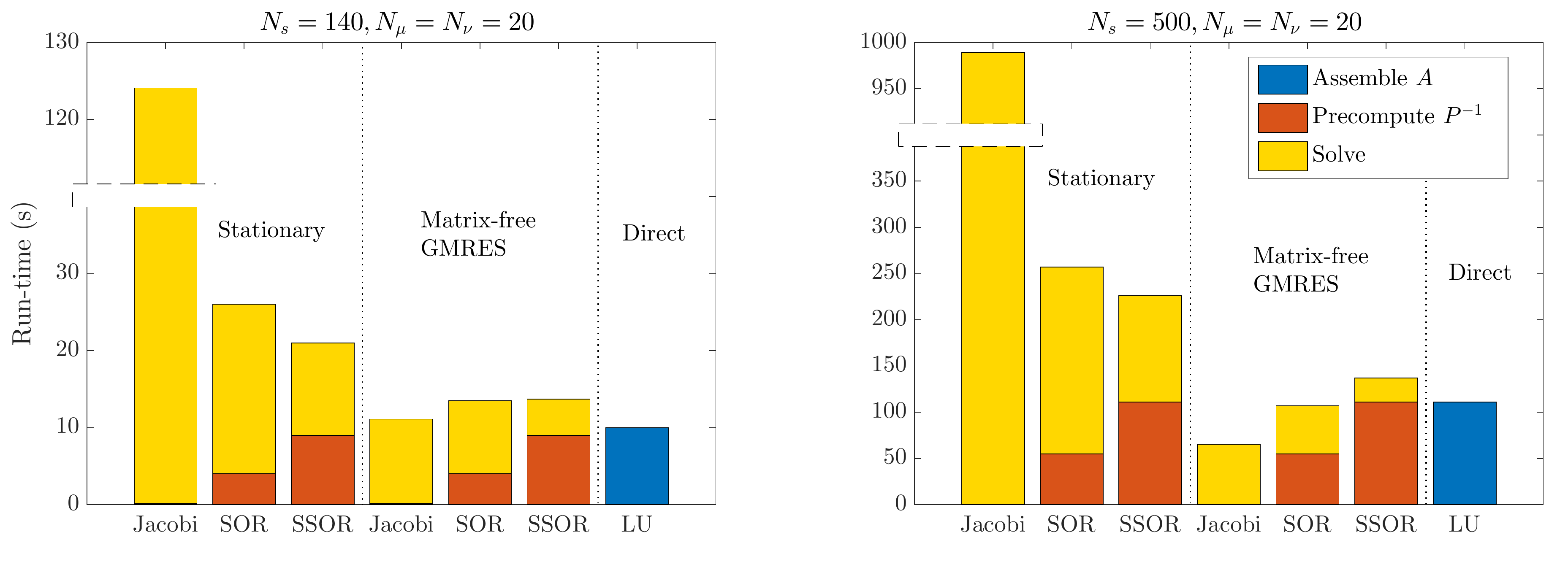}
    \caption{Total run times (in seconds), subdivided in various stages, for the solution of \eqref{crd_linear_system}
    with a DELO-linear formal, $N_\mu=N_\nu=20$, $N_s=140$ (left), and $N_s=500$ (right), according to Tables~\ref{tab:run_times_assemble}-\ref{tab:run_times_mf2}.
    The first three bars correspond to various preconditioned Richardson iterations. The following three bars represent run times for the matrix-free GMRES method, again with various preconditioners. Finally, the last bar represents a direct approach, i.e., an LU solver. Only in this case $A$ is explicitly assembled using Algorithm~\ref{algo} and the ``Solve'' time can hardly be seen being negligible.}
    \label{fig:bar_plot}
\end{figure*}

\begin{figure*}[ht!]
    \centering
    \includegraphics[width=0.95\textwidth]{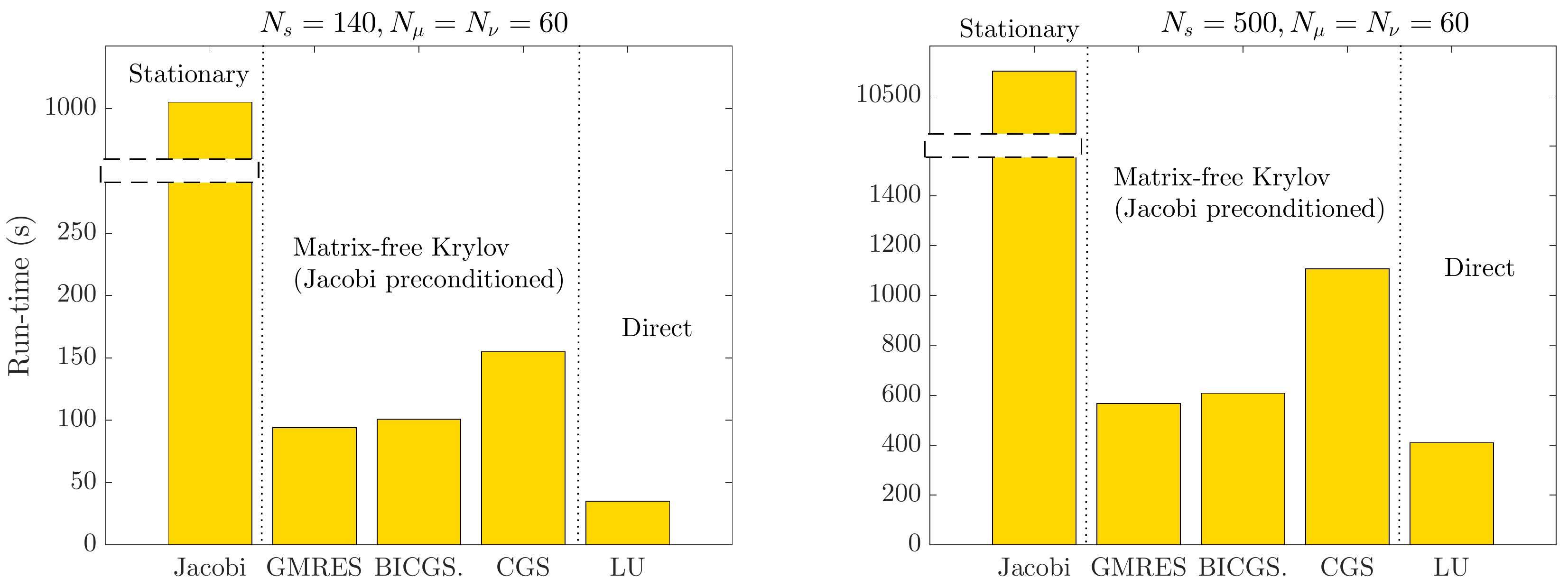}
    \caption{Total run times (in seconds) for the solution of \eqref{crd_linear_system} with a DELO-linear formal, $N_\mu=N_\nu=60$, $N_s=140$ (left), and $N_s=500$ (right). For both plots the first bar corresponds to a stationary Jacobi method. In following three bars the Jacobi method is accelerated with various Krylov methods. Finally, the last bar represents a direct approach, i.e, an LU solver. Only in this case $A$ is explicitly assembled using Algorithm~\ref{algo}.}
    \label{fig:bar_plot2}
\end{figure*}

For run time measures, we also considered the matrix-free application of the iterative methods under investigation.
In this regard, the GMRES method preconditioned with Jacobi has proven to be the most advantageous
choice in terms of convergence and run time (approximately 10 to 20 time faster with respect to standard Jacobi-Richardson, cf. 
Figures~\ref{fig:bar_plot}--\ref{fig:bar_plot2}). 
Since GMRES requires only one matrix-vector multiplication in each iteration, it becomes more advantageous 
when $A$ is less sparse or in the matrix-free case, where the evaluation of $A\xx$ is computationally expensive. 
We remark that parallelism can be extensively exploited both in Algorithm~\ref{algo} and in the application of the Jacobi preconditioner.

We remark that the fully algebraic formulation
of the linear radiative transfer problem
allows us to explicitly assemble the matrix $A$ and to apply direct methods to solve the corresponding linear system.
In this case, the solving time is negligible with respect to to the assembly time. This approach can be more advantageous than the matrix-free one
if high accuracy is required, if \eqref{crd_linear_system} has to be solved for many right hand sides, or if the number of discrete frequency 
and directions is particularly large.
On the other hand, the matrix-free GMRES method preconditioned with Jacobi seems to be preferable if the number of spatial nodes is particularly 
large (e.g., for 3D geometries). 

We finally recall that many available numerical libraries implement matrix-free Krylov routines, enabling an almost painless transition from 
stationary iterative methods to Krylov methods. In terms of implementation, Krylov methods can simply wrap already existing radiative transfer routines 
based on stationary iterative methods.
The application of Krylov methods
to more complex linear radiative transfer problems
will be the subject of study of forthcoming papers.

\begin{table*}[!ht]
\begin{tabular}{c|r|r|r|r}
      $N_\mu=N_\nu$ & $N_s = 40$ & $N_s = 80$ & $N_s = 140$ & $N_s = 500$ \\
      \hline 
      \rule{0pt}{1.1\normalbaselineskip}20 & 1.0 & 3.4 & 10.0 & 111 \\
      30 & 1.7 & 5.4 & 14.7 & 173\\
      40 & 2.6 & 8.0 & 20.5 & 241\\
      50 & 3.6 & 10.4 & 26.8 & 321\\
      60 & 4.8 & 13.1 & 34.7 & 410
      
\end{tabular}
    \caption{Run times (in seconds) to assemble the matrix $A$ in \eqref{crd_linear_system} for the DELO-linear formal solver,
    varying $N_\mu=N_\nu$ and $N_s$.}
    \label{tab:run_times_assemble}
    
    \end{table*}

\begin{table*}[!ht]
    \centering
\parbox{.49\textwidth}{

\begin{tabular}{l|r|r|r|r}
      & Richard. & GMRES & BICGS. & CGS \\ 
       \hline 
\rule{0pt}{1.1\normalbaselineskip}no prec. & - &  33.7 [134] &  71.0 [140] &  71.1 [143]\\
Jacobi & 124 [505] & 10.9 [41] & 13.0 [24] & 19.3 [37] \\
SOR & 21.9 [90] & 9.5 [37] & 8.4 [16] & 17.7 [35] \\
SSOR & 11.9 [51] & 4.7 [17] & 5.5 [10] & 6.3 [12]
\end{tabular}
    \caption{Run times (in seconds) and number of iterations in brackets 
    for the matrix-free solution of \eqref{crd_linear_system}
    for the DELO-linear formal solver,
    with $N_s=140$ and $N_\mu=N_\nu=20$.
    \label{tab:run_times_mf}}



\begin{tabular}{l|r|r|r|r}
      & Richard. & GMRES & BICGS. & CGS \\ 
       \hline 
\rule{0pt}{1.1\normalbaselineskip}no prec. & - & 192 [231] & 304 [171] & 258 [150]  \\
Jacobi & 989 [1391] & 64.7 [71] & 69.7 [38] & 125 [71] \\
SOR &  202 [238]  &  52.0 [59] &  51.4 [28] & 204 [109]\\
SSOR & 115 [132] & 25.8 [28] & 29.4 [16] & 29.0 [16]
\end{tabular}
    \caption{Same as Table~\ref{tab:run_times_mf}, but with $N_s=500$.
    }
    \label{tab:run_times_mf2}
}
\hfill
\parbox{.49\textwidth}{
\vspace*{0.55cm}
      

\begin{tabular}{c|r|r|r|r}
      
      $N_\mu=N_\nu$ & $N_s = 40$ & $N_s = 80$ & $N_s = 140$ & $N_s = 500$ \\
      \hline 
      \rule{0pt}{1.1\normalbaselineskip}20 & 1.5 [18] & 4.4 [29] & 10.6 [41] & 66.0 [71] \\
      30 & 3.5 [18] & 9.5 [29] & 23.1 [41] &  149 [71]  \\
      40 & 6.0 [19] & 16.9 [29] & 39.3 [41] & 253 [71] \\
      50 & 9.1 [19] & 25.0 [29] & 59.5 [41] & 379 [71] \\
      60 & 12.9 [19] & 36.2 [29] & 93.8 [41] & 567 [71] \\
      
\end{tabular}
    \caption{Run times (in seconds) and number of iterations in brackets
    for the matrix-free solution of \eqref{crd_linear_system} using GMRES method
    with the Jacobi preconditioner and varying $N_\mu=N_\nu$ and $N_s$.}
    \label{tab:run_times_mf3}}
\end{table*}


\begin{acknowledgements}
Special thanks are extended to E. Alsina Ballester, N. Guerreiro, and T. Simpson for particularly enriching comments on this work.
The authors gratefully acknowledge the Swiss National Science Foundation (SNSF) for financial support through grant CRSII5$\_180238$. Rolf Krause acknowledges the funding from the European High-Performance Computing Joint Undertaking (JU) under grant agreement No 955701 (project TIME-X). The JU receives support from the European Union's Horizon 2020 research and innovation programme and Belgium, France, Germany, Switzerland. 
\end{acknowledgements}

\FloatBarrier
\twocolumn

\bibliographystyle{aa}
\bibliography{bibfile}

\begin{thebibliography}{25}
\expandafter\ifx\csname natexlab\endcsname\relax\def\natexlab#1{#1}\fi

\bibitem[{Anusha {et~al.}(2009)Anusha, Nagendra, Paletou, \&
  L{\'e}ger}]{anusha2009preconditioned}
Anusha, L., Nagendra, K., Paletou, F., \& L{\'e}ger, L. 2009, The Astrophysical
  Journal, 704, 661

\bibitem[{{Anusha} \& {Nagendra}(2011)}]{anusha2011c}
{Anusha}, L.~S. \& {Nagendra}, K.~N. 2011, \apj, 738, 116

\bibitem[{{Anusha} \& {Nagendra}(2012)}]{anusha2012e}
{Anusha}, L.~S. \& {Nagendra}, K.~N. 2012, \apj, 746, 84

\bibitem[{{Anusha} \& {Nagendra}(2013)}]{anusha2013}
{Anusha}, L.~S. \& {Nagendra}, K.~N. 2013, \apj, 767, 108

\bibitem[{{Anusha} {et~al.}(2011){Anusha}, {Nagendra}, \&
  {Paletou}}]{anusha2011b}
{Anusha}, L.~S., {Nagendra}, K.~N., \& {Paletou}, F. 2011, \apj, 726, 96

\bibitem[{Badri {et~al.}(2019)Badri, Jolivet, Rousseau, \&
  Favennec}]{badri2019preconditioned}
Badri, M., Jolivet, P., Rousseau, B., \& Favennec, Y. 2019, Computers \&
  Mathematics with Applications, 77, 1453

\bibitem[{Castro \& Trelles(2015)}]{castro2015spatial}
Castro, R.~O. \& Trelles, J.~P. 2015, Journal of Quantitative Spectroscopy and
  Radiative Transfer, 157, 81

\bibitem[{Gutknecht(1993)}]{gutknecht1993variants}
Gutknecht, M.~H. 1993, SIAM journal on scientific computing, 14, 1020

\bibitem[{Hubeny \& Burrows(2007)}]{hubeny2007new}
Hubeny, I. \& Burrows, A. 2007, The Astrophysical Journal, 659, 1458

\bibitem[{Ipsen \& Meyer(1998)}]{ipsen1998idea}
Ipsen, I.~C. \& Meyer, C.~D. 1998, The American mathematical monthly, 105, 889

\bibitem[{{Janett} {et~al.}(2021){Janett}, {Benedusi}, {Belluzzi}, \&
  {Krause}}]{janett2021a}
{Janett}, G., {Benedusi}, P., {Belluzzi}, L., \& {Krause}, R. 2021, \aap,
  accepted

\bibitem[{{Janett} {et~al.}(2017{\natexlab{a}}){Janett}, {Carlin}, {Steiner},
  \& {Belluzzi}}]{janett2017a}
{Janett}, G., {Carlin}, E.~S., {Steiner}, O., \& {Belluzzi}, L.
  2017{\natexlab{a}}, \apj, 840, 107

\bibitem[{{Janett} {et~al.}(2017{\natexlab{b}}){Janett}, {Steiner}, \&
  {Belluzzi}}]{janett2017b}
{Janett}, G., {Steiner}, O., \& {Belluzzi}, L. 2017{\natexlab{b}}, \apj, 845,
  104

\bibitem[{{Janett} {et~al.}(2018){Janett}, {Steiner}, \&
  {Belluzzi}}]{janett2018b}
{Janett}, G., {Steiner}, O., \& {Belluzzi}, L. 2018, \apj, 865, 16

\bibitem[{{Klein} {et~al.}(1989){Klein}, {Castor}, {Greenbaum}, {Taylor}, \&
  {Dykema}}]{klein1989}
{Klein}, R.~I., {Castor}, J.~I., {Greenbaum}, A., {Taylor}, D., \& {Dykema},
  P.~G. 1989, \jqsrt, 41, 199

\bibitem[{Lambert {et~al.}(2015)Lambert, Josselin, Ryde, \&
  Faure}]{lambert2015new}
Lambert, J., Josselin, E., Ryde, N., \& Faure, A. 2015, Astronomy \&
  Astrophysics, 580, A50

\bibitem[{{Landi Degl'Innocenti} \&
  {Landolfi}(2004)}]{landi_deglinnocenti+landolfi2004}
{Landi Degl'Innocenti}, E. \& {Landolfi}, M. 2004, Astrophysics and Space
  Science Library, Vol. 307, {Polarization in Spectral Lines} (Dordrecht:
  Kluwer Academic Publishers)

\bibitem[{Meurant \& Duintjer~Tebbens(2020)}]{meurant2020}
Meurant, G. \& Duintjer~Tebbens, J. 2020, Krylov Methods for Nonsymmetric
  Linear Systems: From Theory to Computations (Springer)

\bibitem[{{Nagendra} {et~al.}(2009){Nagendra}, {Anusha}, \&
  {Sampoorna}}]{nagendra2009}
{Nagendra}, K.~N., {Anusha}, L.~S., \& {Sampoorna}, M. 2009, \memsai, 80, 678

\bibitem[{Paletou \& Anterrieu(2009)}]{paletou2009conjugate}
Paletou, F. \& Anterrieu, E. 2009, Astronomy \& Astrophysics, 507, 1815

\bibitem[{Saad(2003)}]{saad2003iterative}
Saad, Y. 2003, Iterative methods for sparse linear systems (SIAM)

\bibitem[{Saad \& Schultz(1986)}]{saad1986gmres}
Saad, Y. \& Schultz, M.~H. 1986, SIAM Journal on scientific and statistical
  computing, 7, 856

\bibitem[{{Sampoorna} {et~al.}(2019){Sampoorna}, {Nagendra}, {Sowmya},
  {Stenflo}, \& {Anusha}}]{sampoorna2019}
{Sampoorna}, M., {Nagendra}, K.~N., {Sowmya}, K., {Stenflo}, J.~O., \&
  {Anusha}, L.~S. 2019, \apj, 883, 188

\bibitem[{Sonneveld(1989)}]{sonneveld1989cgs}
Sonneveld, P. 1989, SIAM journal on scientific and statistical computing, 10,
  36

\bibitem[{Van~der Vorst(1992)}]{van1992bi}
Van~der Vorst, H.~A. 1992, SIAM Journal on scientific and Statistical
  Computing, 13, 631

\end{thebibliography}

\end{document}